\documentclass{amsart}
\usepackage{amsmath,amsfonts,latexsym,amssymb}

\newtheorem{theorem}{Theorem}
\newtheorem{lemma}{Lemma}

\newtheorem{proposition}{Proposition}

\theoremstyle{remark}

\begin{document}

\title[Effective determination of Maass forms]{On effective determination of Maass forms from\\ central values of Rankin-Selberg $L$-function}
\author{Ritabrata Munshi \and Jyoti Sengupta}
\address{School of Mathematics, Tata institute of Fundamental Research, 1 Dr. Homi Bhabha Road, Colaba, Mumbai 400005, India.} 
\email{rmunshi@math.tifr.res.in \and sengupta@math.tifr.res.in}


\begin{abstract}
We address the problem of identifying a Hecke-Maass cusp form $f$ of full level from the central values of the Rankin-Selberg $L$-functions $L(1/2,f\otimes h)$ where $h$ runs through the set of Hecke-Maass eigenforms of full level. We prove a quantitative result in this direction. 
\end{abstract}

\maketitle

%======================================================================================================================
%======================================================================================================================

\section{Introduction}
\label{intro}

The history of determining modular forms from central values of the $L$-function of its twists is fairly long (see \cite{CD}, \cite{GHS}, \cite{Li}, \cite{Liu}, \cite{Luo}, \cite{LR}, \cite{M}, \cite{P} and \cite{Z}).  It was first considered by Luo and Ramakrishnan \cite{LR}.  They showed that if two cuspidal normalised newforms $f$ and $g$ of weight $2k$ (resp. $2k^{'}$) and level $N$ (resp. $N^{'}$), have the property that $L(\frac{1}{2}, f\otimes \chi_d) = L(\frac{1}{2},g\otimes \chi_d)$ for all quadratic characters $\chi_d$ then $k=k^{'}, \ N=N^{'}\ \mbox{and} \ f=g$.  Chinta and Diaconu \cite{CD} generalised this result to self-dual forms on $GL(3)$.  The  next step was to consider $GL(2)$ twists.  In this direction Luo proved the following in \cite{Luo}. Let $f$ and $g$ be as above.  Suppose there exists a  positive integer $\ell$ and infinitely many primes $p$,  such that for all forms $h$ in the Hecke basis $H_{2\ell} (\Gamma_0(p))$ of newforms of weight $2\ell$ and level $p$,
$$
L(\tfrac{1}{2}, f\otimes h)=L(\tfrac{1}{2}, g\otimes h).
$$
Then  $k=k^{'}, \ N=N^{'}\ \mbox{and} \ f=g$. The second author with Ganguly and Hoffstein \cite{GHS} (see also \cite{Z}), extended this result to the weight aspect, namely they considered two newforms $f$ and $g$ of level one and weights $2k \ \mbox{and} \ 2k^{'}$.  If $L(\frac{1}{2}, f\otimes h) = L(\frac{1}{2}, g\otimes h)$ for all normalised cuspidal Hecke eigenforms $h$ of weight $2\ell$  and level one for infinitely many $\ell$, then $k =k^{'}$ and $f=g$.\\

In another direction one may ask if the result of Luo and Ramakrishnan alluded to above can be made effective. This was accomplished  by the first author in \cite{M}.  He showed that if $f,g$ are normalised cuspidal newforms of weight $2k$ (resp. $2k^{'}$) and level $N$ (resp. $N^{'}$) with $k, k^{'} \leq K$ and $N, N^{'}\leq Q$, then for any $\delta > 0 $ and $C$ we have a constant $\gamma(C,\delta)$ such that if $ L(\frac{1}{2}, f\otimes \chi_d) = C \ L(\frac{1}{2}, g\otimes \chi_d)$ for all $d$ in the range $\mid d \mid \leq \gamma(C,\delta)Q^{3/2} (QK)^{4+\delta}$ then $f=g$.\\

The aim of the present article is to obtain an analogue of the above result with $f$ and $g$ being Maass cusp forms for the full modular group, $SL_2(\mathbb Z)$, with the twists being made by Maass cusp forms of varying Laplace eigenvalue. To this end we will prove the following

\begin{theorem}
\label{mthm}
Let $f$, $g$ be normalized Hecke-Maass forms of full level with Laplace eigenvalues $\frac{1}{4}+\mu^2$, $\frac{1}{4}+\nu^2$ respectively, with $|\mu|, |\nu| \leq \Lambda$. For any $\delta> 0$ we have a constant $\gamma(\delta)$ such that if
\begin{align}
\label{cond}
L(1/2,f\otimes h)=L(1/2,g\otimes h),
\end{align} 
for all Hecke-Maass forms $h$ of full level and eigenvalues $\frac{1}{4}+t^2$ with $|t|\leq \gamma(\delta)\Lambda^{4\theta+3+\delta}$ then $f=g$. Here $\theta$ is the exponent towards the Ramanujan conjecture for $GL(2)$ Maass forms of full level.  
\end{theorem}

On the one hand this may be construed as a generalisation of the result in \cite{GHS} to the case of nonholomorphic forms with the weight $k$ being replaced by the eigenvalue of the Laplacian and on the other hand due to its effectivity it can be viewed as a generalisation of the result \cite{M}. We emphasize that in our case the forms $f$ and $g$ are Maass cusp forms.  The case where $f$ and $g$ are both Maass forms and the twisting is done by holomorphic forms has been considered earlier in \cite{P}. However the result there is not effective in our sense. \\

We conclude this section with the note that our result can be generalized to arbitrary level with some extra work. However there are no new technical difficulties. Also we note that by the well-known work of Kim and Sarnak $\theta=\frac{7}{64}$ is admissible. The Ramanujan conjecture predicts that $\theta=0$. Finally, the bound in our result can be improved by using subconvexity results in Section \ref{pop1}, and Voronoi summation formula in Section \ref{pop23}.

%=====================================================================================================================

\section{Preliminaries}
\label{sor}

Let $f:\mathbb H\rightarrow \mathbb C$ be a normalized (i.e. $\lambda_f(1)=1$) Hecke-Maass cusp form (of weight $0$) for $SL(2,\mathbb Z)$, with Laplace eigenvalue $\frac{1}{4}+\mu^2\geq 0$, and with Fourier expansion
$$
f(x+iy)=\sqrt{y}\sum_{n\neq 0}\lambda_f(n)K_{i\mu}(2\pi|n|y)e(nx).
$$     
Here $e(x)=e^{2\pi ix}$ and $K_{i\mu}$ is the MacDonald-Bessel function. The spectral parameter $\mu$ in this case of full level is known to be real $\mu\in\mathbb R$ (so that $\frac{1}{4}+\mu^2\geq \frac{1}{4}$). The associated $L$-function, which for $\sigma>1$ is given by the absolutely convergent Dirichlet series and Euler product
$$
L(s,f)=\sum_{n=1}^\infty \lambda_f(n)n^{-s}=\prod_p \left(1-\alpha_{f,1}(p)p^{-s}\right)^{-1}\left(1-\alpha_{f,2}(p)p^{-s}\right)^{-1},
$$
extends to an entire function, and satisfies a functional equation 
$$
\Lambda(s,f):=\gamma(s,f)L(s,f)=\varepsilon(f)\Lambda(1-s,f).
$$
The gamma factor is given by
$
\gamma(s,f)=\pi^{-s}\Gamma\left(\frac{s+\delta+i\mu}{2}\right)\Gamma\left(\frac{s+\delta-i\mu}{2}\right),
$
where $\delta=0,1$ is the parity of the form. The numbers $\alpha_{f,i}(p)$ are called the local parameters of the form $f$ at the prime $p$. For the sake of exposition we shall only deal with the even Maass cusp forms (i.e. $\lambda_f(-n)=\lambda_f(n)$).\\

Let $u_j$, $j=1,2,\dots$ be an orthogonal basis of the space of even Maass cusp forms, consisting of Hecke-Maass cusp forms. For each $u_j$ we will denote the corresponding normalized Fourier coefficients by $\lambda_j(n)$ and the respective Laplace eigenvalues by $\frac{1}{4}+t_j^2$. Let $h$ be an even function satisfying the following two conditions (i)~$h(t)$ is holomorphic in the strip $|\text{Im}(t)|\leq \frac{1}{2}+\varepsilon$, (ii) $h(t)\ll (1+|t|)^{-2-\varepsilon}$ in the same strip. Let 
\begin{align*}
h_0=&\frac{2}{\pi}\int_0^\infty th(t)(\tanh \pi t)\:\mathrm{d}t\\
h^+(x)=&2i\int_{-\infty}^\infty tJ_{2it}(x)h(t)(\cosh \pi t)^{-1}\: \mathrm{d}t\\
h^-(x)=&\frac{4}{\pi}\int_0^\infty t K_{2it}(x)h(t)(\sinh \pi t)\: \mathrm{d}t.
\end{align*}
Then we have the following trace formula of Kuznetsov: 
\begin{lemma}
\label{ktf}
For two positive integers $n$ and $m$, we have
\begin{align*}
\sum_{j=1}^\infty h(t_j)\omega_j\lambda_j(n)\lambda_j(m)=&\delta_{n,m}h_0+\sum_{c=1}^\infty\frac{1}{2c}\left\{S(m,n;c)h^+\left(\tfrac{4\pi\sqrt{mn}}{c}\right)+S(-m,n;c)h^-\left(\tfrac{4\pi\sqrt{mn}}{c}\right)\right\}\\
&+\frac{1}{4\pi}\int_{\mathbb R}h(t)\omega(t)\left(\frac{m}{n}\right)^{it}\sum_{m_0|m}m_0^{-2it}\sum_{n_0|n}n_0^{2it}\:\mathrm{d}t,
\end{align*} 
where $\delta_{m,n}$ is the Kronecker symbol, 
$$
S(m,n;c)=\sideset{}{^\star}\sum_{\alpha\bmod{c}}e\left(\frac{\alpha m+\bar\alpha n}{c}\right)
$$ 
is the classical Kloosterman sum,
\begin{align*}
\omega_j=4\pi\|u_j\|^{-2}(\cosh \pi t_j)^{-1},
\end{align*}
and
\begin{align*}
\omega(t)=4\pi^{2+2it}\Gamma\left(\tfrac{1}{2}+it\right)^{-2}\zeta(1+2it)^{-2}(\cosh \pi t)^{-1}.
\end{align*}\\
\end{lemma}

Let $f$ be as above. For $h$ another even Hecke-Maass cusp form for $SL(2,\mathbb Z)$, with Laplace eigenvalue $\frac{1}{4}+t^2$, and local parameters $\alpha_{h,i}(p)$, we define the Rankin-Selberg $L$-function $L(s,f\otimes h)$ by the absolutely convergent Euler product 
$$
L(s,f\otimes h)=\prod_p \prod_{i,j=1,2}\left(1-\alpha_{f,i}(p)\alpha_{h,j}(p)p^{-s}\right)^{-1}
$$
in the right half-plane $\text{Re}(s)>1$. In this half plane it is also given by an absolutely convergent Dirichlet series 
$$
L(s,f\otimes h)=\zeta(2s)\sum_{n=1}^\infty \lambda_f(n)\lambda_h(n)n^{-s}.
$$
It is well-know that this $L$-function extends to a meromorphic function, with a pole at $s=1$ only in the case $f=\bar h$ ($=h$ for full level). Otherwise the function is entire. Moreover we have the functional equation 
$$
\Lambda(s,f\otimes h):=\gamma(s;\mu,t)L(s,f\otimes h)=\Lambda(1-s,f\otimes h),
$$
where the gamma factor is given by
$$
\gamma(s;\mu,t)=\pi^{-2s}\Gamma\left(\tfrac{s+i(\mu+t)}{2}\right)\Gamma\left(\tfrac{s+i(\mu-t)}{2}\right)\Gamma\left(\tfrac{s-i(\mu-t)}{2}\right)\Gamma\left(\tfrac{s-i(\mu+t)}{2}\right).
$$
The gamma factor depends only on the spectral parameters $\mu$ and $t$. A consequence of the functional equation is the following expression for the central value, called the approximate functional equation.
\begin{lemma}
For $f$ and $h$ as above, we have
\begin{align}
\label{afe}
L(\tfrac{1}{2},f\otimes h)=2\sum_{n=1}^\infty \frac{\lambda_f(n)\lambda_h(n)}{\sqrt{n}}V_{\mu,t}\left(n\right)
\end{align}
where
\begin{align}
\label{v-fun}
V_{\mu,t}(y)=\frac{1}{2\pi i}\int_{(3)}y^{-u}\left(\cos \frac{\pi u}{4A}\right)^{-16A}\zeta(1+2u)\frac{\gamma(\tfrac{1}{2}+u;\mu,t)}{\gamma(\tfrac{1}{2};\mu,t)}\frac{\mathrm{d}u}{u},
\end{align}
and $A$ is any positive integer. \\
\end{lemma}

The cosine function in \eqref{v-fun} gives an exponential decay in the integral as $|\text{Im}(u)|\rightarrow\infty$. The size of the decay can be regulated by choosing $A$ large enough. Suppose $\mu$ and $t$ are such that $|\mu|,|t| \leq T$. Then using Stirling approximation for the gamma functions and moving the contour to the right we get that $V(n)\ll_N T^{-N}$ for any $N>0$ if $n\gg T^{2+\varepsilon}$. This reflects the well-known fact that approximate functional equation has effective length given by the square root of the size of the conductor, which in our case is $\ll T^4$.\\

Let 
\begin{align}
\label{wt-fun}
\mathcal H_{T,M}(t)=e^{-\frac{(t-T)^2}{M^2}}+e^{-\frac{(t+T)^2}{M^2}}
\end{align}
where the positive parameters $T$ and $M\leq T^{1-\varepsilon}$ will be chosen conveniently later. Observe that the modified weight function $\mathcal H_{T,M}(t)V_{\mu,t}(n)$, as function in $t$ satisfies the properties required for the Kuznetsov trace formula (see Lemma~\ref{ktf}). Also note that the weight function $\mathcal H_{T,M}(t)$ localizes $t$ at size $T$. Suppose $|\mu|\leq \Lambda\ll T^{1-\varepsilon}$ and $|t|\asymp T$. Moving the contour in \eqref{v-fun} to $\varepsilon$, and applying the Stirling approximation 
\begin{align}
\label{stirling}
\Gamma(z)=e^{(z-\frac{1}{2})\log z-z+\frac{1}{2}\log 2\pi}\left(1+\frac{1}{12z}+\frac{1}{288z^2}+O\left(\frac{1}{|z|^3}\right)\right)
\end{align}
to the gamma functions appearing in \eqref{v-fun}, we get that 
\begin{align}
\label{v-fun-2}
V_{\mu,t}(y)=V_{\mu,t}^\flat(y)+O\left(\frac{T^{\varepsilon}}{T^3}\right).
\end{align}
The new function $V_{\mu,t}^\flat(y)$ has the additional `scaling property' that 
\begin{align}
\label{scaling-prop}
\frac{\partial^j}{\partial \tau^j}V_{\mu,T+\tau M}^\flat(y)\ll_j T^{\varepsilon}
\end{align}
in the range $|\tau|\ll T^{\varepsilon}$.\\

For smaller values of $y$ we can evaluate $V_{\mu,t}^\flat(y)$ (or $V_{\mu,t}(y)$) precisely by shifting the contour in \eqref{v-fun} to the left upto $\text{Re}(u)=-\frac{1}{2}+\varepsilon$. We pass through a double pole at $u=0$, and the residue there gives us the leading term. The integral over the line $\text{Re}=-\frac{1}{2}+\varepsilon$ yields that error term. Using the first term approximation from \eqref{stirling}, we get
\begin{align}
\label{v-fun-22}
V_{\mu,t}^\flat(y)=\log\left(\frac{|t|}{2\sqrt{y}}\right)+\gamma-1+O(y^{\frac{1}{2}+\varepsilon}|t|^{-1}),
\end{align}
where $\gamma$ is the Euler constant. Also the implied constant is absolute. The same asymptotic holds for $V_{\mu,t}(y)$. \\

%Next we shall recall the Voronoi summation formula for $SL(2,\mathbb Z)$ automorphic forms. Let $f$ be a Maass cusp form as above %with Laplace eigenvalue $\frac{1}{4}+\mu^2\geq 0$. We will use the following Voronoi type summation formula (see Meurman %\cite{Me}).  
%\begin{lemma}
%Let $h$ be compactly supported smooth function on $(0,\infty)$. We have
%\begin{align}
%\label{voronoi2}
%\sum_{n=1}^\infty \lambda_f(n)e_q\left(an\right)\Psi(n)=\frac{1}{q}\sum_{\pm}\sum_{n=1}^\infty \lambda_f(\mp %n)e_q\left(\pm\bar{a}n\right)\Psi^{\star,\pm}\left(\frac{n}{q^2}\right)
%\end{align}
%where $\bar{a}$ is the multiplicative inverse of $a\bmod{q}$, and
%\begin{align*}
%\Psi^{\star,-}(y)=&-\frac{\pi}{\cosh \pi\mu}\int_0^\infty \Psi(x)\{Y_{2i\mu}+Y_{-2i\mu}\}\left(4\pi\sqrt{xy}\right)dx\\
%\Psi^{\star,+}(y)=&4\cosh \pi\mu\int_0^\infty \Psi(x)K_{2i\mu}\left(4\pi\sqrt{xy}\right)dx.
%\end{align*}
%\end{lemma}

%====================================================================================================================
%====================================================================================================================

\section{Proof of Theorem \ref{mthm}}
\label{pot}

To prove the main theorem we will relate the the Fourier coefficients of the forms $f$, $g$ by computing the twisted average of the central values of the Rankin-Selberg $L$-functions. More precisely we will prove the following:
\begin{proposition}
\label{p2}
Let $f$ be a normalized even Hecke-Maass form of full level and with spectral parameter $\mu$ (so that the Laplace eigenvalue is $\frac{1}{2}+\mu^2$), and $|\mu|\leq \Lambda$. Suppose $\{u_j\}_j$ is a orthogonal basis of the space of even Maass forms of full level. Let $p$ be a prime. Let $\mathcal H_{T,M}(t)$ be as defined in \eqref{wt-fun}, with $T$ such that $\Lambda^{1+\varepsilon}\ll T$ and $M=T^{1-\varepsilon}$. Let $\omega_j$ be as in Lemma \ref{ktf}. Then we have
\begin{align}
\label{eqp2}
\sum_{j=1}^\infty \mathcal H_{T,M}(t_j)\omega_jL(\tfrac{1}{2},f\otimes u_j)\lambda_j(p)
=\frac{4M\lambda_f(p)}{\pi\sqrt{p}}\mathop{\int}_{\mathbb R} &e^{-t^2}(T+tM)\left[\log_0\left(\frac{T+tM}{2\sqrt{p}}\right)+\gamma-1\right] \mathrm{d}t\\
\nonumber &+O\left(\sqrt{p}\:\Lambda\:T^{1+\varepsilon}\right)
\end{align} 
where the implied constant depends only on $\varepsilon$. (Here $\log_0(x)=\log x$ for $x>1$ and $=0$ otherwise.)
\end{proposition}
  
The proof of the above proposition is rather delicate and forms the technical heart of this paper. In another vein we will show that if $f\neq g$ then they cannot have same Fourier coefficients. In fact the distinction starts quite early in the sequence. The following proposition gives a quantitative result in this direction.

\begin{proposition}
\label{p1}
Let $f$, $g$ be Hecke-Maass cusp forms of weight $0$, full level and Laplace eigenvalues $\frac{1}{4}+\mu^2$ and $\frac{1}{4}+\nu^2$ respectively. Suppose $|\mu|, |\nu| \leq \Lambda$ and $f\neq g$. Then there is a prime $p\ll \Lambda^{2+\varepsilon}$, such that
\begin{align*}
|\lambda_f(p)- \lambda_g(p)|\gg \Lambda^{-4\theta-\varepsilon},
\end{align*} 
where $\theta$ is the exponent towards the Ramanujan conjecture for Maass forms. \\
\end{proposition}

We will now prove the main theorem assuming the two propositions. Let $f$ and $g$ be two Hecke-Maass forms of full level and spectral parameters $\mu$ and $\nu$ respectively. Suppose $|\mu|,|\nu|\leq \Lambda$. Choose $T\gg \Lambda^{1+\varepsilon}$,  $M=T^{1-\varepsilon}$ and pick any prime $p\ll \Lambda^{2+\varepsilon}$. (Recall our convention that the value $\varepsilon$ may be different at each occurrence.) Suppose 
$$
L(\tfrac{1}{2},f\otimes u_j)=L(\tfrac{1}{2},g\otimes u_j)
$$
for $u_j$ as above satisfying $|t_j|\ll T^{1+\varepsilon}$. Then we get 
$$
\sum_{j=1}^\infty \mathcal H_{T,M}(t_j)\omega_jL(\tfrac{1}{2},f\otimes u_j)\lambda_j(p)=\sum_{j=1}^\infty \mathcal H_{T,M}(t_j)\omega_jL(\tfrac{1}{2},g\otimes u_j)\lambda_j(p)+O(T^{-N}).
$$ 
Applying Proposition \ref{p2} it now follows that
$$
|\lambda_f(p)-\lambda_g(p)|\frac{4M}{\pi\sqrt{p}}\mathop{\int}_{\mathbb R} e^{-t^2}(T+tM)\left[\log_0\left(\frac{T+tM}{2\sqrt{p}}\right)+\gamma-1\right] \mathrm{d}t= O\left(\sqrt{p}\:\Lambda\:T^{1+\varepsilon}\right)
$$
and consequently
$$
|\lambda_f(p)-\lambda_g(p)|= O\left(\frac{p\Lambda}{T}\:T^{\varepsilon}\right)= O\left(\frac{\Lambda^3}{T}\:T^{\varepsilon}\right).
$$
This contradicts Proposition \ref{p1} if $T\gg \Lambda^{4\theta+3+\varepsilon}$. This concludes the proof of Theorem~\ref{mthm}.  

%==================================================================================================================
%==================================================================================================================
%==================================================================================================================

\section{Distance between Fourier coefficients of distinct forms}
\label{pop1}

In this section we will prove Proposition \ref{p1}. Let $f$ and $g$ be as in the statement of the proposition. Let us write 
\begin{align*}
\lambda_f(n)-\lambda_g(n)=b(n),
\end{align*}
and suppose that $|b(p)|\leq B$ for all primes $p\leq P$. Then we claim the following for any square-free integer 
$n \leq P$:
\begin{align*}
|b(n)| \leq 2\tau(n)^2|n|^{\theta}B,
\end{align*}
where $\tau(n)$ is the divisor function and $\theta$ is the exponent towards the Ramanujan conjecture for Maass forms. So we are assuming that $|\lambda_h(n)|\leq \tau(n)|n|^{\theta}$ for any Hecke-Maass form $h$ of full level. By the work of Kim and Sarnak we know that $\theta= \frac{7}{64}$ is admissible. \\

The bound on $b(n)$ can be obtained by induction on the number of prime factors in $n$. Suppose $n=p_1p_2\dots p_k$. If $k=1$ then the statement is clear. Otherwise let $m=p_1\dots p_{k-1}$ and using multiplicativity of the Fourier coefficient we obtain
\begin{align*}
|\lambda_f(n)-\lambda_g(n)| &\leq  |\lambda_f(m)||\lambda_f(p_k)-\lambda_g(p_k)|
+|\lambda_f(m)-\lambda_g(m)||\lambda_g(p_k)|.
\end{align*}
Now using the induction hypothesis and bound for individual Fourier coefficients, we get that the above expression is
\begin{align*} 
\leq  \tau(m)|m|^{\theta}B+2\tau(m)^2|m|^{\theta}B2|p_k|^\theta\leq 2\tau(n)^2|n|^\theta B.
\end{align*} 
This proves our assertion. 
\\

To derive a consequence of this bound, let $F(y)$ be a non-negative smooth bounded function supported in $[\frac{1}{2},1]$. Then for any $X\leq P$, we get
\begin{align}
\label{rs}
\left|\sideset{}{^\flat}\sum_{n}\lambda_f(n)^2F\Bigl(\frac{n}{X}\Bigr)
-\sideset{}{^\flat}\sum_{n}\lambda_g(n)\lambda_f(n)F\Bigl(\frac{n}{X}\Bigr)\right|
&=\left|\sideset{}{^\flat}\sum_{n}b(n)\lambda_f(n)F\Bigl(\frac{n}{X}\Bigr)\right|\\
\nonumber \leq 2B\sum_{n}\tau(n)^3|n|^{2\theta}F\Bigl(\frac{n}{X}\Bigr)&\ll BX^{1+2\theta}(\log X)^4,
\end{align}
where the implied constant depends only on $F$. Here the $\flat$ denotes that the sum is restricted to square-free integers. \\

Next we will employ the Rankin-Selberg $L$-function to estimate the sums appearing on the left hand side of (\ref{rs}). To this end consider the Dirichlet series
\begin{align}
\label{s}
D_{h,f}(s)=\sideset{}{^\flat}\sum_{}\lambda_h(n)\lambda_f(n)n^{-s}
=\prod_p\Bigl(1+\frac{\lambda_h(p)\lambda_f(p)}{p^s}\Bigr),
\end{align}
where $h$ is either $f$ or $g$. We can relate this Dirichlet series with the Rankin-Selberg $L$-function $L(s,h\otimes f)$. Indeed expanding the local Euler factors we have
$$
L(s,h\otimes f)=\prod_p \prod_{i,j}\left(1-\alpha_{h,i}(p)\alpha_{f,j}(p)p^{-s}\right)^{-1}=\prod_p\left(1+\lambda_{h}(p)\lambda_{f}(p)p^{-s}+\dots\right),
$$
where the local factors of $D_{h,f}(s)$ appear as the first degree approximation. Using the bound $|\alpha_{h,i}(p)|\leq p^{\theta}$ (where $\theta$ is as above) for the local parameters , we conclude that we have a Dirichlet series $L_{h,f}(s)$ which converges absolutely in the region $\sigma\geq \frac{1}{2}+2\theta+\varepsilon$, and in this region it satisfies the bound $L_{h,f}(s)\ll_{\varepsilon} 1$, and such that
\begin{align*}
D_{h,f}(s)=L(s,h\otimes f)L_{h,f}(s).
\end{align*}\\

For any compactly supported smooth function $F$, it follows by Mellin inversion and contour shifting (Perron's formula) that 
\begin{align}
\label{ss}
\sideset{}{^\flat}\sum_{n}\lambda_h(n)\lambda_f(n)F\Bigl(\frac{n}{X}\Bigr)
&=c_FL_{h,f}(1)\mathop{\text{Res}}_{s=1}\:L(s,h\otimes f)X+O(\Lambda^{1-4\theta}X^{\frac{1}{2}+2\theta+\varepsilon}),
\end{align}
where the constant $c_F=\int F(y)dy$, and the implied constant depends only on $F$. The error term comes from the integral over the contour given by $\sigma=\frac{1}{2}+2\theta+\varepsilon$ in the $s$-plane. The integral converges due to the rapid decay of the Mellin transform of $F$ (as $F$ is compactly supported in $[1/2,1]$ and smooth). The Rankin-Selberg $L$-function is bounded by $$L(\tfrac{1}{2}+2\theta+\varepsilon+it,h\otimes f)\ll (\Lambda(3+|t|))^{1-4\theta},$$ which is the well-known convexity bound. \\

Comparing (\ref{rs}) and (\ref{ss}) we obtain the following bound for any $X\leq P$, 
\begin{align*}
c_FL_{f,f}(1)\mathop{\text{Res}}_{s=1}\:L(s,f\otimes f)=O(\Lambda^{1-4\theta}X^{-\frac{1}{2}+2\theta+\varepsilon}+BX^{2\theta}(\log X)^4).
\end{align*}
Now clearly $L_{f,f}(1)\gg 1$ (as the Dirichlet series defining $L_{f,f}(s)$ converges absolutely for $\sigma\geq \frac{1}{2}+2\theta+\varepsilon$, and $\theta$ can be taken to be $\frac{7}{64}$), and consequently the left hand side is bounded below by $\Lambda^{-\varepsilon}$. So we are led to a contradiction if $X\asymp \Lambda^{2+\varepsilon}$, and $B\asymp \Lambda^{-4\theta-\varepsilon}$. This completes the proof of Proposition~\ref{p2}.

%==================================================================================================================
%==================================================================================================================
%==================================================================================================================

\section{Twisted first moment: The main term}
\label{pop2}

The rest of the paper is devoted to proving Proposition \ref{p2}. Using approximate functional equation \eqref{afe}, we get
$$
\Xi_{T,M}(f,p):=\sum_{j=1}^\infty \mathcal H_{T,M}(t_j)\omega_jL(\tfrac{1}{2},f\otimes u_j)\lambda_j(p)=2\sum_{n=1}^\infty \frac{\lambda_f(n)}{\sqrt{n}}\sum_{j=1}^\infty \mathcal H_{T,M}(t_j)V_{\mu,t_j}(n)\omega_j\lambda_j(n)\lambda_j(p),
$$
where the weight function $\mathcal H_{T,M}(t)$ is defined in \eqref{wt-fun}. Using the decomposition \eqref{v-fun-2}, and the rapid decay of both the functions $\mathcal H_{T,M}(t)V_{\mu,t}(n)$ and $\mathcal H_{T,M}(t)V_{\mu,t}^\flat(n)$ for $n\gg T^{2+\varepsilon}$, one obtains 
$$
\Xi_{T,M}(f,p)=2\sum_{n=1}^\infty \frac{\lambda_f(n)}{\sqrt{n}}\sum_{j=1}^\infty \mathcal H_{T,M}(t_j)V_{\mu,t_j}^\flat(n)\omega_j\lambda_j(n)\lambda_j(p)+O\left(\frac{MT^{\varepsilon}}{T}\right),
$$
where the implied constant depends only on $\varepsilon$.\\

Set $h_n(t)=\mathcal H_{T,M}(t)V_{\mu,t}^\flat(n)$. (For notational simplicity we are suppressing the parameters $\mu$, $T$ and $M$.) Next we apply Kuznetsov trace formula (Lemma \ref{ktf}) to get
\begin{align}
\label{kuz-split}
\Xi_{T,M}(f,p)=\mathcal D +\mathcal E +\mathcal J+\mathcal K,
\end{align}
where 
$$
\mathcal D=\frac{4}{\pi}\frac{\lambda_f(p)}{\sqrt{p}}\int_0^\infty th_n(t)(\tanh \pi t) \mathrm{d}t
$$
is the diagonal term, 
$$
\mathcal E=2\sum_{n=1}^\infty \frac{\lambda_f(n)}{\sqrt{n}}\frac{1}{4\pi}\int_{\mathbb R}h_n(t)\omega(t)\left(\frac{p}{n}\right)^{it}\sum_{m_0|p}m_0^{-2it}\sum_{n_0|n}n_0^{2it}\:\mathrm{d}t
$$
is the contribution of the continuous spectrum, and 
\begin{align*}
\mathcal J&=\sum_{n=1}^\infty \frac{\lambda_f(n)}{\sqrt{n}}\sum_{c=1}^\infty \frac{S(n,p;c)}{c}h^+\left(\frac{4\pi\sqrt{np}}{c}\right)\\
\mathcal K&=\sum_{n=1}^\infty \frac{\lambda_f(n)}{\sqrt{n}}\sum_{c=1}^\infty \frac{S(-n,p;c)}{c}h^-\left(\frac{4\pi\sqrt{np}}{c}\right)
\end{align*}
are the terms related to the $J$-Bessel function and the $K$-Bessel function respectively.\\

The diagonal $\mathcal D$ will yield the main term. The remaining terms will contribute to the error. Applying \eqref{v-fun-22}, and using the definition of the weight function $\mathcal H_{T,M}(t)$ we get
$$
\mathcal D=\frac{4\lambda_f(p)}{\pi\sqrt{p}}\int_0^\infty e^{-\frac{(t-T)^2}{M^2}}t\left[\log\left(\frac{t}{2\sqrt{p}}\right)+\gamma-1\right] \mathrm{d}t+O(p^{\theta-\frac{1}{2}+\varepsilon}T).
$$
Again $\theta$ is the exponent towards the Ramanujan conjecture. Performing a change of variables we get our next result. \\

\begin{lemma}
\label{diagonal-cont}
We have
$$
\mathcal D=\frac{4M\lambda_f(p)}{\pi\sqrt{p}}\mathop{\int}_{\mathbb R} e^{-t^2}(T+tM)\left[\log_0\left(\frac{T+tM}{2\sqrt{p}}\right)+\gamma-1\right] \mathrm{d}t+O(p^{\theta-\frac{1}{2}+\varepsilon}T)
$$
where the implied constant depends only on $\theta$. (Recall that $\log_0 x=\log x$ for $x>1$ and $=0$ elsewhere.)
\end{lemma}

Hence the diagonal term is of size $\asymp MT\log T$. Also note that the leading term does not depend on the spectral parameter $\mu$, and the implied constant in the error term is uniform with respect to this parameter.
%==================================================================================================================
%==================================================================================================================
%==================================================================================================================

\section{Contribution from the continuous spectrum}
\label{pop2con}

Next we will estimate the Eisenstein series (continuous spectrum) contribution
\begin{align*}
\mathcal E=&2\sum_{n=1}^\infty \frac{\lambda_f(n)}{\sqrt{n}}\frac{1}{4\pi}\int_{\mathbb R}\mathcal H_{T,M}(t)V_{\mu,t}^\flat(n)\omega(t)\left(\frac{p}{n}\right)^{it}\sum_{m_0|p}m_0^{-2it}\sum_{n_0|n}n_0^{2it}\:\mathrm{d}t.
\end{align*} 
The Dirichlet series
$$
\sum_{n=1}^\infty\frac{\lambda_f(n)}{n^{\frac{1}{2}+u+it}}\sum_{n_0|n}n_0^{2it}
$$
is associated with the Rankin-Selberg convolution of $f$ with an Eisenstein series. More precisely (in the region of absolute convergence) we have
$$
\sum_{n=1}^\infty\frac{\lambda_f(n)}{n^{\frac{1}{2}+u+it}}\sum_{n_0|n}n_0^{2it}=\frac{L\left(\tfrac{1}{2}+u+it,f\right)L\left(\tfrac{1}{2}+u-it,f\right)}{\zeta(1+2u)}.
$$
So using \eqref{v-fun-2} it follows that 
\begin{align*}
&\sum_{n=1}^\infty\frac{\lambda_f(n)}{n^{\frac{1}{2}+it}}V_{\mu,t}^\flat(n)\sum_{n_0|n}n_0^{2it}=\sum_{n=1}^\infty\frac{\lambda_f(n)}{n^{\frac{1}{2}+it}}V_{\mu,t}(n)\sum_{n_0|n}n_0^{2it}+O\left(\frac{T^{\varepsilon}}{T^2}\right)\\
=&\frac{1}{2\pi i}\int_{(3)}\left(\cos \frac{\pi u}{4A}\right)^{-16A}\frac{\gamma(\tfrac{1}{2}+u;\mu,t)}{\gamma(\tfrac{1}{2};\mu,t)}L\left(\tfrac{1}{2}+u+it,f\right)L\left(\tfrac{1}{2}+u-it,f\right)\frac{\mathrm{d}u}{u}+O\left(\frac{T^{\varepsilon}}{T^2}\right).
\end{align*}\\

The analytic conductor of the $L$-function $L\left(\tfrac{1}{2}+u\pm it,f\right)$ is bounded by $3+(|u|+|t|+|\mu|)^2$. We move the contour to $\text{Re}(u)=\varepsilon$, and consider the integral over $t$, which is given by
\begin{align*}
&\int_{\mathbb R}\mathcal H_{T,M}(t)\omega(t)\frac{\gamma(\tfrac{1}{2}+u;\mu,t)}{\gamma(\tfrac{1}{2};\mu,t)}L\left(\tfrac{1}{2}+u+it,f\right)L\left(\tfrac{1}{2}+u-it,f\right)p^{it}\sum_{m_0|p}m_0^{-2it}\:\mathrm{d}t.
\end{align*}
Applying Cauchy, Stirling's approximation and the definitions of the functions $\mathcal H_{T,M}(t)$ and $\omega(t)$, the job of estimating the above integral reduces to getting bounds for 
\begin{align*}
\mathop{\int}_{\frac{T}{10}}^{10T}\left|L\left(\tfrac{1}{2}+u+it,f\right)\right|^2\:\mathrm{d}t.
\end{align*}
Applying approximate functional equation we see that upto a negligible error term the $L$-value is given by a partial series of length $|u|+T$. (Recall that $|\mu|\leq \Lambda \leq T^{1-\varepsilon}$.) This leads us to consider 
\begin{align}
\label{prep-ls}
\mathop{\int}_{\frac{T}{10}}^{10T}\left|\sum_{n\leq |u|+T}\frac{\lambda_f(n)}{n^{\frac{1}{2}+u+it}}W_{\mu,t}(n)\right|^2\:\mathrm{d}t,
\end{align}
where the smooth weight function is given by
$
W_{\mu,t}(y)=\frac{1}{2\pi i}\int_{(3)}y^{-v}\mathcal W_{\mu,t}(v)\mathrm{d}v,
$
with
\begin{align*}
\mathcal W_{\mu,t}(v)=\left(\cos \frac{\pi v}{4A}\right)^{-16A}\frac{\Gamma\left(\frac{\frac{1}{2}+u+v+i(t+\mu)}{2}\right)\Gamma\left(\frac{\frac{1}{2}+u+v+i(t-\mu)}{2}\right)}{\Gamma\left(\frac{\frac{1}{2}+u+i(t+\mu)}{2}\right)\Gamma\left(\frac{\frac{1}{2}+u+i(t-\mu)}{2}\right)}\frac{1}{v}.
\end{align*} 
Moving the contour to $\text{Re}(v)=\varepsilon$ and interchanging the order of summation and integration, we see that \eqref{prep-ls} is dominated by
\begin{align*}
\mathop{\int}_{\frac{T}{10}}^{10T}\left|\int_{(\varepsilon)}\mathcal W_{\mu,t}(v)\sum_{n\leq |u|+T}\frac{\lambda_f(n)}{n^{\frac{1}{2}+u+v+it}}\mathrm{d}v\right|^2\:\mathrm{d}t.
\end{align*}
Applying Cauchy (to the inner integral) we bound this by
\begin{align*}
\mathop{\int}_{\frac{T}{10}}^{10T}\int_{(\varepsilon)}|\mathcal W_{\mu,t}(v_1)||\mathrm{d}v_1|\int_{(\varepsilon)}|\mathcal W_{\mu,t}(v_2)|\left|\sum_{n\leq |u|+T}\frac{\lambda_f(n)}{n^{\frac{1}{2}+u+v_2+it}}\right|^2|\mathrm{d}v_2|\:\mathrm{d}t.
\end{align*}
Then applying Stirling (to obtain a point-wise bound for $|\mathcal W_{\mu,t}(v)|$), and the large sieve inequality (see \cite{IK})
$$
\int_{T}^{2T}\left|\sum_{n\leq N}a_nn^{it}\right|^2\ll (T+N)\sum_{n\leq N}|a_n|^2,
$$
together with $\sum_{n\leq N}|\lambda_f(n)|^2\ll |\mu|^{\varepsilon}N$, we get that \eqref{prep-ls} is dominated by $(|u|+T)T^{\varepsilon}$. We summarize the outcome of our analysis in the following\\

\begin{lemma}
\label{eisenstein-cont}
Let $\mathcal E$ be as given in \eqref{kuz-split}. We have
$$
\mathcal E\ll T(pT)^{\varepsilon},
$$
where the implied constant depends only on $\varepsilon$.
\end{lemma}

%==================================================================================================================
%==================================================================================================================
%==================================================================================================================

\section{The $J$-Bessel contribution}
\label{pop22}

Next we turn our attention to the analysis of the contribution of the sums involving Kloosterman sums. First consider the part with $J$-Bessel function, i.e. the term $\mathcal J$ of \eqref{kuz-split}. We will closely follow the method of Sarnak \cite{S}, and its modifications by Lau-Liu-Ye \cite{LLY} and Li \cite{L}. Recall that $h_n(t)=\mathcal H_{T,M}(t)V_{\mu,t}^\flat(n)$ and consider the integral transform
\begin{align*}
h_n^+\left(\frac{4\pi\sqrt{np}}{c}\right)=&2i\int_{-\infty}^\infty tJ_{2it}\left(\frac{4\pi\sqrt{np}}{c}\right)h_n(t)(\cosh \pi t)^{-1}\: \mathrm{d}t.
\end{align*}
We have already noted that $h_n(t)$, and hence $h_n^+\left(\frac{4\pi\sqrt{np}}{c}\right)$, is negligibly small for $n\gg T^{2+\varepsilon}$. So we focus our attention on the complimentary range. In this case 
$$
\frac{\sqrt{np}}{c}\ll \sqrt{p}\:T^{1+\varepsilon}. 
$$ 
We choose $M$ so that $\sqrt{p}\ll MT^{-\varepsilon}$. In the rest of the section we write $x$ in place of $\frac{4\pi\sqrt{np}}{c}$, and show that the integral $h_n^+(x)$ is negligibly small for $|x|\ll T^{1-\varepsilon}M$. \\ 

%Moving the line of integration to $\text{Im}(t)=-\frac{1}{2}+\varepsilon$ and using the bound (for $\tau\in\mathbb R$)
%$$
%J_{2i\tau+1-\varepsilon}(x)\ll \left(\frac{x}{1+|\tau|}\right)^{1-\varepsilon}e^{\pi|\tau|}
%$$
%one gets
%\begin{align*}
%h_n^+(x)=&2i\int_{-\infty}^\infty tJ_{2it}(x)h_n(t)(\cosh \pi t)^{-1}\: \mathrm{d}t.
%\end{align*}\\
Splitting the integral in the definition of $h_n^+(x)$, and using the identity (see \cite{GR} 8.411-11)
$$
\frac{J_{2it}(x)-J_{-2it}(x)}{\cosh \pi t}=-\frac{2i}{\pi}\tanh \pi t\mathop{\int}_{-\infty}^{\infty}\cos(x\cosh \zeta)e\left(\frac{t\zeta}{\pi}\right)\mathrm{d}\zeta
$$
we get
\begin{align}
\label{h+fun}
h_n^+(x)=&\frac{4}{\pi}\mathop{\int}_{0}^\infty th_n(t)\tanh \pi t\mathop{\int}_{-\infty}^{\infty}\cos(x\cosh \zeta)e\left(\frac{t\zeta}{\pi}\right)\mathrm{d}\zeta\: \mathrm{d}t.
\end{align}
Using the definition of $h_n(t)$ and applying integration-by-parts with respect to $\zeta$ once, it follows that
\begin{align*}
h_n^+(x)=&\frac{4}{\pi}\mathop{\int}_{0}^\infty te^{-\frac{(t-T)^2}{M^2}}V_{\mu,t}^\flat(n)\tanh \pi t\mathop{\int}_{-T^{\varepsilon}}^{T^{\varepsilon}}\cos(x\cosh \zeta)e\left(\frac{t\zeta}{\pi}\right)\mathrm{d}\zeta\: \mathrm{d}t +O_N(T^{-N})
\end{align*}
for any $N>0$. By changing variables $\frac{t-T}{M}\rightarrow t$ we get
\begin{align*}
h_n^+(x)=\frac{4M}{\pi}\mathop{\int}_{-\frac{T}{M}}^\infty &(T+tM)e^{-t^2}V_{\mu,T+tM}^\flat(n)\tanh \pi (T+tM)\\
&\times \mathop{\int}_{-T^{\varepsilon}}^{T^{\varepsilon}}\cos(x\cosh \zeta)e\left(\frac{(T+tM)\zeta}{\pi}\right)\mathrm{d}\zeta\: \mathrm{d}t +O_N(T^{-N}).
\end{align*}
At a cost of negligible error term we can now extend the integral over $t$ to $(-\infty,\infty)$. Then we write $h_n^+(x)=h^+_{n,1}(x)+h^+_{n,2}(x)+O(T^{-N})$ where
\begin{align*}
h^+_{n,1}(x)=\frac{4MT}{\pi}\mathop{\int}_{-\infty}^\infty &e^{-t^2}V_{\mu,T+tM}^\flat(n)\tanh \pi (T+tM)\mathop{\int}_{-T^{\varepsilon}}^{T^{\varepsilon}}\cos(x\cosh \zeta)e\left(\frac{(T+tM)\zeta}{\pi}\right)\mathrm{d}\zeta\: \mathrm{d}t.
\end{align*}
We will only treat the integral $h^+_{n,1}(x)$. The other integral $h^+_{n,2}(x)$ can be handled in a similar fashion and in fact can be shown to be of smaller order of magnitude. \\

Let
\begin{align}
\label{k-fun}
k(t)=e^{-t^2}V_{\mu,T+tM}^\flat(n)\tanh \pi (T+tM),
\end{align}
and suppose $\hat k(y)=\int k(t)e(-ty)\mathrm{d}t$ denote the Fourier transform. From the properties of $V_{\mu,T+tM}^\flat(n)$ that we noted in Section \ref{sor} (in particular \eqref{scaling-prop}), it follows that $k^{(j)}(t)\ll_j 1$ and hence $\hat k(y)$ is essentially supported in $[-T^{\varepsilon},T^{\varepsilon}]$ (i.e. the function is negligibly small outside this support). Also it follows that $\hat k^{(j)}(y)\ll_j T^{\varepsilon}$. We have
\begin{align*}
h^+_{n,1}(x)=\frac{4MT}{\pi}\mathop{\int}_{-T^{\varepsilon}}^{T^{\varepsilon}}\hat k\left(-\frac{M\zeta}{\pi}\right)\cos(x\cosh \zeta)e\left(\frac{T\zeta}{\pi}\right)\mathrm{d}\zeta.
\end{align*}
After a change of variables we get
\begin{align*}
h^+_{n,1}(x)=4T\mathop{\int}_{-\frac{MT^{\varepsilon}}{\pi}}^{\frac{MT^{\varepsilon}}{\pi}}\hat k\left(\zeta\right)\cos\left(x\cosh \frac{\pi\zeta}{M}\right)e\left(-\frac{T\zeta}{M}\right)\mathrm{d}\zeta.
\end{align*}\\

Now we can again extend the integral to $(-\infty,\infty)$ at the cost of a negligible error term. Set
\begin{align*}
h^\star_{n,1}(x)=2T\mathop{\int}_{-\infty}^{\infty}\hat k\left(\zeta\right)e\left(\frac{x}{2\pi}\cosh \frac{\pi\zeta}{M}-\frac{T\zeta}{M}\right)\mathrm{d}\zeta,
\end{align*}
so that $h^+_{n,1}(x)=h^\star_{n,1}(x)+h^\star_{n,1}(-x)+O(T^{-N})$. Let 
$$
\phi(\zeta)=\frac{x}{2\pi}\cosh \frac{\pi\zeta}{M}-\frac{T\zeta}{M}.
$$
Then 
\begin{align*}
\phi'(\zeta)&=\frac{x}{2M}\sinh \frac{\pi\zeta}{M}-\frac{T}{M}\\
&=-\frac{T}{M}+\frac{\pi x}{2M^2}\zeta +\text{smaller order terms}.
\end{align*}
Hence if $|x|\leq T^{1-\varepsilon}M$ we get that $\phi'(\zeta)\leq -\frac{T}{10M}$. In particular $\phi(\zeta)$ is monotonic for these values of $x$. Also taking higher order derivatives we get that $|\phi^{(j)}(\zeta)|\ll \frac{|x|}{M^j}\ll \frac{T}{T^{\varepsilon}M^{j-1}}$ for $j\geq 2$. So for $|x|\leq T^{1-\varepsilon}M$, applying the change of variables $\phi(\zeta)\rightarrow \xi$ and then repeated integration-by-parts it follows that $h^\star_{n,1}(x)$ is negligibly small. Thus we arrive at our next result.
\\

\begin{lemma}
\label{J-Bessel-cont}
Let $\mathcal J$ be as given in \eqref{kuz-split}. We have
$$
\mathcal J\ll_N T^{-N}
$$
for any $N>0$.
\end{lemma}

%==================================================================================================================
%==================================================================================================================
%==================================================================================================================

\section{The $K$-Bessel contribution}
\label{pop23}

To estimate the contribution coming from the $K$-Bessel function, i.e. the term $\mathcal K$ in \eqref{kuz-split}, we will follow a similar recipe as in the previous section. Consider 
\begin{align*}
h^-_n(x)=\frac{4}{\pi}\int_{0}^\infty tK_{2it}(x)h_n(t)\sinh \pi t\: \mathrm{d}t.
\end{align*}
Using the integral representation (see \cite{GR} 8.432-4)
$$
K_{2it}(x)=\frac{1}{2}(\cosh \pi t)^{-1}\mathop{\int}_{-\infty}^{\infty}\cos(x\sinh \zeta)e\left(-\frac{t\zeta}{\pi}\right)\mathrm{d}\zeta
$$
we get
\begin{align*}
h^-_n(x)=&\frac{2}{\pi}\mathop{\int}_{0}^\infty th_n(t)\tanh \pi t\mathop{\int}_{-\infty}^{\infty}\cos(x\sinh \zeta)e\left(-\frac{t\zeta}{\pi}\right)\mathrm{d}\zeta\: \mathrm{d}t.
\end{align*}
This can be compared with the integral representation \eqref{h+fun} for $h^+_n(x)$. Next we follow exactly the same steps to arrive at the decomposition $h^-_n(x)=h^-_{n,1}(x)+h^-_{n,2}(x)+O(T^{-N})$ where
\begin{align*}
h^-_{n,1}(x)=\frac{2MT}{\pi}\mathop{\int}_{-\infty}^\infty &e^{-t^2}V_{\mu,T+tM}^\flat(n)\tanh \pi (T+tM)\mathop{\int}_{-T^{\varepsilon}}^{T^{\varepsilon}}\cos(x\sinh \zeta)e\left(-\frac{(T+tM)\zeta}{\pi}\right)\mathrm{d}\zeta\: \mathrm{d}t.
\end{align*}
As before we shall only treat the integral $h^-_{n,1}(x)$. The other integral $h^-_{n,2}(x)$ can be handled in a similar fashion and can be shown to be of smaller order of magnitude. \\

Now we have
\begin{align*}
h^-_{n,1}(x)=\frac{2MT}{\pi}\mathop{\int}_{-T^{\varepsilon}}^{T^{\varepsilon}}\hat k\left(\frac{M\zeta}{\pi}\right)\cos(x\sinh \zeta)e\left(-\frac{T\zeta}{\pi}\right)\mathrm{d}\zeta,
\end{align*}
where $k$ is as defined in \eqref{k-fun}. It again follows that $h^-_{n,1}(x)=h^\sharp_{n,1}(x)+h^\sharp_{n,1}(-x)+O(T^{-N})$ with
\begin{align*}
h^\sharp_{n,1}(x)=T\mathop{\int}_{-\infty}^{\infty}\hat k\left(\zeta\right)e\left(\frac{x}{2\pi}\sinh \frac{\pi\zeta}{M}-\frac{T\zeta}{M}\right)\mathrm{d}\zeta.
\end{align*}
Let 
$$
\psi(\zeta)=\frac{x}{2\pi}\sinh \frac{\pi\zeta}{M}-\frac{T\zeta}{M}.
$$
Then 
\begin{align*}
\psi'(\zeta)&=\frac{x}{2M}\cosh \frac{\pi\zeta}{M}-\frac{T}{M}\\
&=-\frac{T}{M}+\frac{x}{2M}+\frac{\pi^2x}{4M^3}\zeta^2+\text{smaller order terms}.
\end{align*}
Hence if $x\notin [10^{-1}T,10T]$ then as in previous section we can show that $h^\sharp_{n,1}(x)$ is negligibly small. Otherwise using the expansion of $\sinh$ and taking first term approximation we get
\begin{align}
\label{h-sharp}
h^\sharp_{n,1}(x)=T\mathop{\int}_{-\infty}^{\infty}\hat k\left(\zeta\right)e\left(\frac{x\zeta}{2M}\right)e\left(-\frac{T\zeta}{M}\right)\mathrm{d}\zeta+O\left(\frac{T|x|}{M^3}\right).
\end{align} 
\\

It remains to analyse the case where $x=\frac{4\pi\sqrt{np}}{c}\asymp T$. In this case we have $c\asymp \frac{\sqrt{np}}{T} \ll \sqrt{p}T^{\varepsilon}$. The sum is given by
$$
\sum_{n=1}^{\infty}\frac{\lambda_f(n)}{\sqrt{n}}\Psi\left(\frac{n}{N}\right)\sum_{c\ll \sqrt{p}\:T^{\varepsilon}}\frac{S(-n,p\:; c)}{c}h^\sharp_{n,1}\left(\frac{4\pi\sqrt{np}}{c}\right),
$$
where $\Psi$ is a smooth non-negative function on $(0,\infty)$ with compact support, and $\frac{T^2}{p}\ll N\ll T^{2+\varepsilon}$.  The error term from \eqref{h-sharp} gives an error term of size 
$$
O\left(\frac{T^{3+\varepsilon}\sqrt{p}}{M^3}\right),
$$
and the leading term yields
$$
T\mathop{\int}_{-\infty}^{\infty}\hat k\left(\zeta\right)e\left(\frac{T\zeta}{M}\right)\sum_{n=1}^{\infty}\frac{\lambda_f(n)}{\sqrt{n}}\Psi\left(\frac{n}{N}\right)\sum_{c\ll \sqrt{p}\:T^{\varepsilon}}\frac{S(-n,p\:;c)}{c}e\left(\frac{2\pi\sqrt{np}\zeta}{cM}\right)\mathrm{d}\zeta.
$$
The Weil bound is not sufficient for our purpose. Interchanging the order of summations and opening the Kloosterman sum we arrive at
$$
T\mathop{\int}_{-\infty}^{\infty}\hat k\left(\zeta\right)e\left(\frac{T\zeta}{M}\right)\sum_{c\ll \sqrt{p}\:T^{\varepsilon}}\frac{1}{c}\;\sideset{}{^\star}\sum_{a\bmod{c}}e\left(\frac{\bar ap}{c}\right)\left[\sum_{n=1}^{\infty}\frac{\lambda_f(n)}{\sqrt{n}}e\left(\frac{an}{c}\right)\Psi\left(\frac{n}{N}\right)e\left(\frac{2\pi\sqrt{np}\zeta}{cM}\right)\right]\mathrm{d}\zeta.
$$
Taking $M=T^{1-\varepsilon}$, $|\mu|\leq \Lambda$, applying partial
summation and the bound (see \cite{I})
$$
\sum_{n\leq x}\lambda_f(n)e(\alpha n)\ll \sqrt{x}\Lambda(x\Lambda)^{\varepsilon}
$$
we get
$$
\sum_{n=1}^{\infty}\frac{\lambda_f(n)}{\sqrt{n}}e\left(\frac{an}{c}\right)\Psi\left(\frac{n}{N}\right)e\left(\frac{2\pi\sqrt{np}\zeta}{cM}\right)\ll \Lambda\left(1+\frac{\sqrt{p}}{c}\right)T^{\varepsilon}.
$$
%Applying Voronoi summation on the sum over $n$ we get
%$$
%\sum_{n=1}^{\infty}\frac{\lambda_f(n)}{\sqrt{n}}e\left(\frac{an}{c}\right)\Psi\left(\frac{n}{N}\right)e\left(\frac{\sqrt{np}\zeta}
%{2cM}\right)=\frac{1}{c}\sum_{\pm}\sum_{n=1}^\infty \lambda_f(\mp n)e\left(\frac{\pm\bar{a}n}{c}\right)\hat %\Psi^{\pm}\left(\frac{n}{c^2}\right)
%$$
%where 
%\begin{align*}
%\hat\Psi^{-}(y)=&-\frac{\pi}{\cosh \pi\mu}\int_0^\infty \Psi\left(\frac{x}{N}\right)e\left(\frac{\sqrt{xp}\zeta}
%{2cM}\right)\{Y_{2i\mu}+Y_{-2i\mu}\}\left(4\pi\sqrt{xy}\right)\frac{dx}{\sqrt{x}}\\
%\hat\Psi^{+}(y)=&4\cosh \pi\mu\int_0^\infty \Psi\left(\frac{x}{N}\right)e\left(\frac{\sqrt{xp}\zeta}
%{2cM}\right)K_{2i\mu}\left(4\pi\sqrt{xy}\right)\frac{dx}{\sqrt{x}}.
%\end{align*}
\\

\begin{lemma}
\label{K-Bessel-cont}
Let $\mathcal K$ be as given in \eqref{kuz-split}. We have
$$
\mathcal K\ll \sqrt{p}\:\Lambda\:T^{1+\varepsilon}
$$
where the implied constant depends only on $\varepsilon$.
\end{lemma}

Proposition \ref{p2} now follows from Lemma \ref{diagonal-cont}, Lemma \ref{eisenstein-cont}, Lemma \ref{J-Bessel-cont} and Lemma \ref{K-Bessel-cont}.

%========================================================================================================================
%========================================================================================================================
%========================================================================================================================
%========================================================================================================================

\end{document}